\documentclass[a4,11pt]{article}
\usepackage{amsmath,amssymb,amstext}
\usepackage{graphicx}

\usepackage[tmargin=1.1in,bmargin=1.1in,rmargin=1.25in,lmargin=1.25in]{geometry}

\title{Sparse Matrix Decompositions and Graph Characterizations}
\author{Kshitij Khare\\ {\small University of Florida} \and Bala Rajaratnam\\ {\small Stanford University}}

\date{}

\begin{document}

\maketitle

\begin{abstract}
The question of when zeros (i.e., sparsity) in a positive definite matrix $A$ are preserved in its Cholesky decomposition, and vice versa, was addressed by Paulsen et al. \cite{plsnpwrsth} [see {\it Journ. of Funct. Anal.}, {\bf 85}, 151-178]. In particular, they prove that for the pattern of zeros in $A$ to be retained in the Cholesky decomposition of $A$,  the pattern of zeros in $A$ has to necessarily correspond to a chordal (or decomposable) graph associated with a 
specific type of vertex ordering. This result therefore yields a characterization of 
chordal graphs in terms of sparse positive definite matrices. It has also proved to be 
extremely useful in probabilistic and statistical analysis of Markov random fields where zeros in positive definite correlation matrices are intimately related to the notion of stochastic independence. Now, consider a positive definite matrix $A$ and its Cholesky decomposition given by $A = LDL^T$, where $L$ is lower triangular with unit diagonal entries, and $D$ a diagonal matrix with positive entries.  In this paper, we prove that a necessary and sufficient condition for zeros (i.e., sparsity) in a positive definite matrix $A$ to be preserved in its associated 
Cholesky matrix $L$, \, and in addition also preserved in the inverse of the Cholesky matrix 
$L^{-1}$, is that the pattern of zeros corresponds to a co-chordal or homogeneous graph associated with a specific type of vertex ordering. We proceed to provide a second characterization of this class of graphs in terms of determinants of submatrices that correspond to cliques in the 
graph. These results add to the growing body of literature in the field of sparse matrix decompositions, and also prove to be critical ingredients in the probabilistic analysis of an important class of Markov random fields.
\end{abstract}



\section{Introduction}

\noindent
Chordal and co-chordal graphs, and their relationships to sparse matrix decompositions, play an important role in the probabilistic and statistical analysis of Markov random fields (see \cite{khrrjtmcvt, khrrjtmcwp, khrrjtmwdg, roveratocs}). In these models the above classes of graphs are used to encode zeros in covariance or correlation matrices (or their inverses). The zero entries in these positive definite correlation matrices are intimately related to the notion of stochastic independence. 

A characterization of chordal graphs or decomposable graphs, the class of graphs containing no induced cycle of length greater than or equal to $4$, in terms of appropriate sub-manifolds of positive definite matrices was provided in \cite{plsnpwrsth}. In particular, positive definite matrices with zero entries according to a decomposable graph necessarily preserve these zero entries in their respective Cholesky matrices. The task undertaken in this paper is to find parallel and useful characterizations of co-chordal or homogeneous graphs, the class of graphs containing no induced $4$-cycle or $4$-path, in terms of appropriate sub-manifolds of positive definite matrices. 

\indent
Let $G = (V,E)$ denote an undirected graph, where $V = \{1,2, \cdots, |V|\}$ represents 
the finite vertex set and $E$ denotes the corresponding edge set. We use the notation 
$\mathbb{M}_p$ to denote the set of $p \times p$ symmetric matrices and 
$\mathbb{M}_p^+$ to denote the set of $p \times p$ positive definite matrices. Without 
loss of generality, the notation used in this paper specifies the permutation or ordering 
$\sigma \in S_p$, where $S_p$ denotes the symmetric group, by a $p$-tuple describing where 
$(1,2, \cdots, p)$ is sent by $\sigma$. Thus, $\sigma = (1 \;2 \; 5 \; 4 \; 3)$ means 
$\sigma(1) = 1, \sigma(2) = 2, \sigma(3) = 5, \sigma(4) = 4$ and $\sigma(5) = 3$. Without 
ambiguity, in some places in the paper we will denote $\sigma$, an element of the symmetric 
group on $p$-letters, by a $p$-tuple describing where $(u,v,w, \cdots)$ is sent by 
$\sigma$. As we explain shortly, these orderings play an important role in our results. 
Given a graph $G = (V,E)$ and an ordering $\sigma$ of the vertices of the graph, we define 
$$
P_{G_\sigma} = \left\{ \Sigma \in \mathbb{M}_{|V|}^+: \; \; \Sigma_{ij} = 0 
\mbox{ whenever } (\sigma^{-1} (i), \sigma^{-1} (j)) \notin E \right\}, 
$$

\noindent
and 
$$
\mathcal{L}_{G_\sigma} = \left\{ L \in \mathbb{M}_{|V|}: \; \; L_{ii} = 1, L_{ij} = 0 
\mbox{ for } i < j \mbox{ or } (\sigma^{-1} (i), \sigma^{-1} (j)) \notin E \right\}. 
$$

\noindent
The space $P_{G_\sigma}$ is essentially a sub-manifold of the space of $|V| \times |V|$ 
positive definite matrices where the elements are restricted to be zero whenever the 
corresponding edge (under the ordering $\sigma$) is missing from $E$. Similarly, the space 
$\mathcal{L}_{G_\sigma}$ is a subspace of lower triangular matrices with diagonal entries 
equal to $1$, such that the elements in the lower triangle are restricted to be zero 
whenever the corresponding edge (under the ordering $\sigma$) is missing from $E$. We now 
state the main theorem of the paper. It characterizes co-chordal or homogeneous graphs in 
terms of (1) sparse matrix decompositions and (2) determinants of submatrices of cliques in the 
graph. 
\newtheorem{thm}{Theorem}
\begin{thm} \label{hmgnschrzn}
Consider a graph $G = (V,E)$ together with an ordering of its vertices as denoted by 
$\sigma$. Then the following statements are equivalent. 
\begin{enumerate}
\item $G$ is a homogeneous graph and $\sigma$ is a Hasse tree based elimination scheme\footnote{a certain type of vertex ordering that will be formally defined later in the paper.}. 
\item If $D$ is an arbitrary diagonal matrix with positive diagonal entries, then 
$$
L \in \mathcal{L}_{G_\sigma} \Leftrightarrow L^{-1} \in \mathcal{L}_{G_\sigma} 
\Leftrightarrow \Sigma := LDL^T \in P_{G_\sigma}. 
$$
\item Let $\Sigma \in P_{G_\sigma}$ be arbitrarily chosen. Let $\Sigma = LDL^T$ denote 
its modified Cholesky decomposition, where $L$ is a lower triangular matrix with unit
diagonal entries and $D$ is a diagonal matrix with diagonal entries $D_{ii}, \; i = 1,2, 
\cdots, p$. Then for any maximal clique $C$ of the graph $G$, 
$$
\left| (\Sigma^{-1})_{\sigma(C)} \right| = \prod_{i \in \sigma(C)} \frac{1}{D_{ii}}. 
$$
\end{enumerate}
\end{thm}

\indent
The outline of the remainder of the paper is as follows. Section \ref{preliminaries} introduces terminology and notation from both linear algebra and graph theory that is required in subsequent sections. Section \ref{first_characterization} provides a first characterization of co-chordal graphs in terms of sparse matrix decompositions.  Section \ref{second_characterization} provides a second characterization of co-chordal graphs in terms of determinants of sub-matrices. The results in Sections \ref{first_characterization} and \ref{second_characterization} are illustrated through examples, which a sophisticated reader can skip.

\section{Preliminaries} \label{preliminaries}

\subsection{Graph theory}

\noindent
This section introduces notation and terminology that is required in 
subsequent sections.
An undirected graph $G=(V,E)$ consists of two sets $V$ and $E$, with $V$ 
representing the set of vertices, and $E\subseteq V \times V$ the set 
of edges satisfying :
$$ \;(u,v)\in E\,\iff \,(v,u)\in E $$

\noindent
When $(u,v)\in E$, we say that $u$ and $v$ are \textit{adjacent} in $G$. A graph is 
said to be \textit{complete} if all the vertices are adjacent to each other, i.e., 
$(u,v) \in E$ for all $u, v \in V$ such that $u \neq v$. A subgraph of $V$ induced by 
$A \subset V$ is the graph $G' = (A, E \cap (A \times A))$. 
\newtheorem{definition}{Definition}
\begin{definition}
A \textit{path} connecting  two  distinct vertices $u$ and $v$ in $G$ is a sequence of 
distinct vertices $\left(u_0,u_1,\ldots,u_n\right)$ where $u_0=u$ and $u_n=v$, and 
for every  $i=0,\ldots,n-1$, $(u_i, u_{i+1}) \in E$. 
\end{definition}

\begin{definition}
A \textit{cycle} is a path with an additional edge between the two endpoints $u_0$ and $u_n$.
\end{definition}

\begin{definition}
A set of vertices $A \subset V$ is said to constitute a \textit{clique} 
if the graph induced by $A$ is a complete subgraph of $V$. Equivalently, 
a clique is a set of vertices in $V$ which are all adjacent to each other. 
\end{definition}

\begin{definition}
A set of vertices $A \subset V$ is said to be a \textit{ maximal clique} 
if $A$ is a clique and is not contained in another clique. Equivalently, 
$A \subset V$ is a  maximal clique if it is a clique and the graph induced 
by $A \cup \{u\}$, for any $u \in V \setminus A$, is no longer a clique. 
\end{definition}

\subsection{Modified Cholesky decomposition} \label{mdfdchldmp}

\noindent
If $\Sigma$ is a positive definite matrix, then there exists a unique decomposition 
\begin{equation} \label{inverselwr}
\Sigma = LDL^T, 
\end{equation}

\noindent
where $L$ is a lower triangular matrix with unit diagonal entries and $D$ a diagonal matrix 
with positive diagonal entries. This decomposition of $\Sigma$ is referred to as the 
{\it modified Cholesky decomposition} of $\Sigma$ (see \cite{pourahmadi}). The lemma below 
provides an explicit formulation of the inverse of a lower triangular matrix with unit 
diagonal entries, and will be useful in subsequent sections. 
\newtheorem{lemma}{Lemma}
\begin{lemma} \label{inverselrt}
Let $L$ be a $p \times p$ lower triangular matrix with diagonal entries equal to $1$. Let 
$$
\mathcal{A} = \cup_{r=2}^p \left\{ {\bf \tau}: {\bf \tau} \in \{1, 2, \cdots, p\}^r, 
\tau_i < \tau_{i-1} \; \; \forall \; 2 \leq i \leq r \right\}, 
$$

\noindent
and 
$$
L_{\bf \tau} = \prod_{i=2}^{dim({\bf \tau})} L_{\tau_{i-1} \tau_i} \;, \; {\bf \tau} \in 
\mathcal{A}, 
$$

\noindent
where $dim(\tau)$ denotes the length of the vector $\tau$. Then $L^{-1} = N$, where 
\begin{eqnarray*}
N_{ij} = \begin{cases}
0 & \mbox{ if } i < j \cr 
1 & \mbox{ if } i = j \cr 
\sum_{{\bf \tau} \in \mathcal{A}, \tau_1 = i, \tau_{dim({\bf \tau})} = j} 
(-1)^{dim({\bf \tau})-1} L_{\bf \tau} & \mbox{ if } i > j. 
\end{cases}
\end{eqnarray*}
\end{lemma}

\subsection{Decomposable graphs} \label{dcmpgphpty}

\begin{figure}
\begin{center}
\includegraphics[width=3in,height=3in]{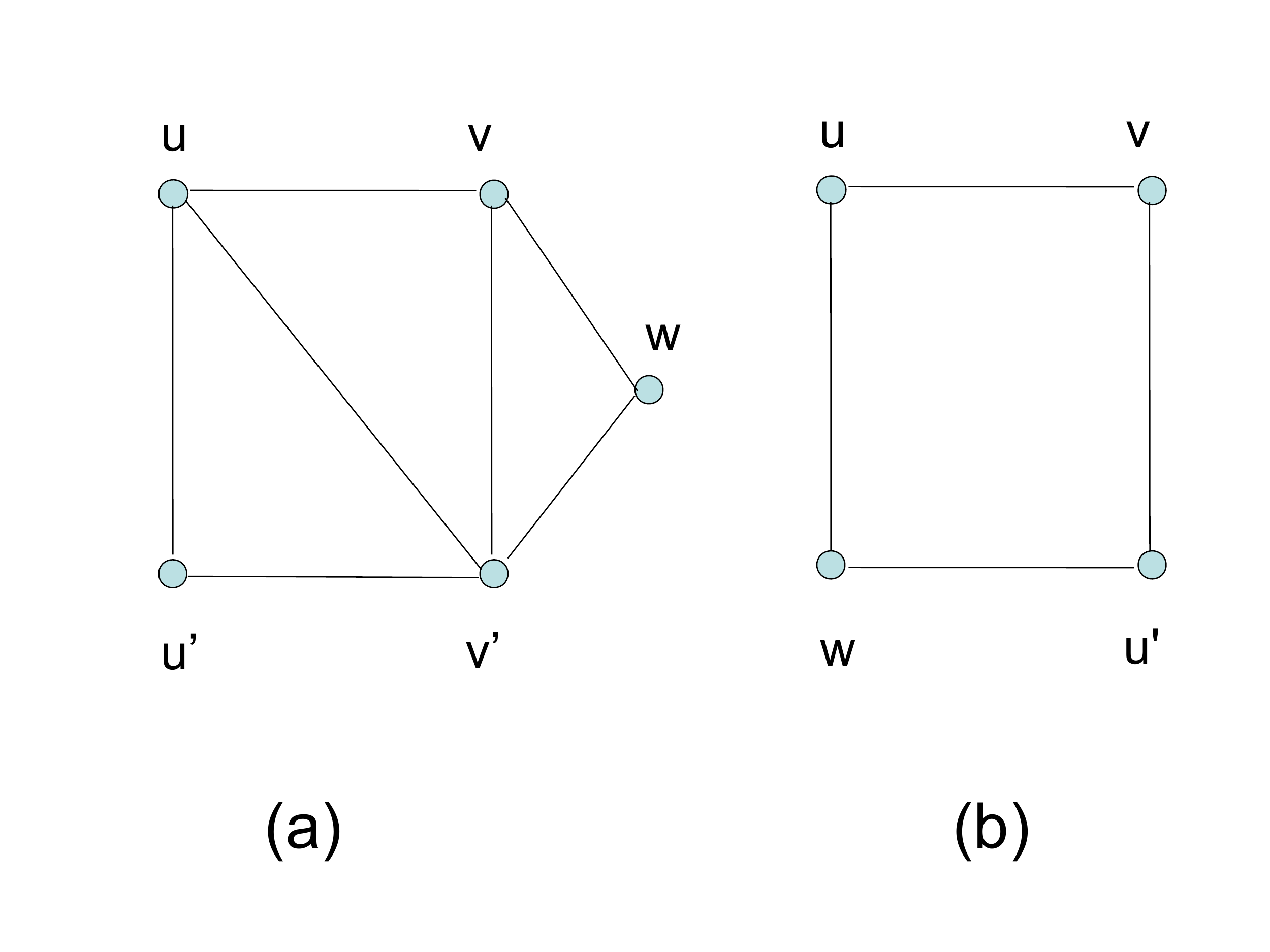}
\end{center}
\caption{(a) A decomposable graph, and (b) A non-decomposable graph}
\label{figure1}
\end{figure}

\noindent
An undirected graph $G$ is said to be {\it decomposable} if any induced subgraph does not 
contain a cycle of length greater than or equal to four. They are also sometimes known as 
chordal graphs or triangulated graphs. See Figure \ref{figure1} for an example of a 
decomposable graph and a non-decomposable graph. Since their introduction by Chvatal 
\cite{dcmplgrchv}, these graphs have been well studied, and are used in various fields such 
as optimization, computer science, probability and statistics. An important branch of 
probability and statistics where the class of decomposable graphs has proven to be quite 
useful is the study of Markov random fields/Graphical models. Decomposable graphs have several 
characterizations. One such characterization is in terms of vertex orderings. We first 
introduce notation and terminology that is required in order to formally state this 
characterization. 

\begin{definition}
For an undirected graph $G = (V,E)$, an ordering $\sigma$ of $V$ is known as a perfect 
vertex elimination scheme for $G$ if for every triplet $i,j,k$ with $1 \leq i < j < k 
\leq p$ the following holds. 
$$
(\sigma^{-1} (j), \sigma^{-1} (i)) \in E, (\sigma^{-1} (k), \sigma^{-1} (i)) \in  E 
\Rightarrow (\sigma^{-1} (k), \sigma^{-1} (j)) \in E. 
$$
\end{definition}

\noindent
A perfect vertex elimination scheme $\sigma$ for the decomposable graph $G$ in Figure 
\ref{figure1} {\it (a)} is given by $\sigma: (u,u',v,v',w) \stackrel{\sigma}{\rightarrow} 
(3,4,2,5,1)$. 

\indent
The existence of such an ordering characterizes decomposable graphs (see 
Paulsen et al \cite{plsnpwrsth}). More formally, an undirected graph $G = (V,E)$ is 
decomposable iff there exists an ordering $\sigma$ of $V$, which is a perfect vertex 
elimination scheme. For a given decomposable graph $G = (V,E)$, there can however be 
several orderings which gives rise to perfect vertex elimination schemes. A constructive 
way to obtain such an ordering is given in Lauritzen \cite{lrtzngphmd}. There is an 
interesting and useful connection between decomposable graphs, orderings which give rise 
to perfect vertex elimination schemes, and the matrix spaces $P_{G_\sigma}$ and 
$\mathcal{L}_{G_\sigma}$. 
\begin{lemma}[Paulsen et al \cite{plsnpwrsth}] \label{dcmplsprpl}
Let $G = (V,E)$ be a decomposable graph, and $\sigma$ an ordering of $V$ which corresponds 
to a perfect vertex elimination scheme for $G$. Then for any positive definite matrix 
$\Sigma$ with modified Cholesky decomposition given by $\Sigma = LDL^T$, the following 
holds. 
$$
L \in \mathcal{L}_{G_\sigma} \Leftrightarrow \Sigma \in P_{G_\sigma}. 
$$
\end{lemma}

\noindent
Hence, for $\Sigma \in P_{G_\sigma}$, the zeros in $\Sigma$ are preserved in the lower 
triangle of the corresponding matrix $L$ obtained from the modified Chloesky decomposition. 
Moreover for $L \in \mathcal{L}_{G_\sigma}$, the zeros in $L$ are preserved in the matrix 
$\Sigma$ obtained by $\Sigma = LDL^T$, for any diagonal matrix $D$ with positive diagonal 
entries. The converse of Lemma \ref{dcmplsprpl} is also true. 
\begin{lemma}[Paulsen et al \cite{plsnpwrsth}] \label{cnvrsdmpgr}
Let $G = (V,E)$ be a graph, $\sigma$ be an ordering of $V$, and $D$ be an arbitrary diagonal 
matrix with positive diagonal entries. Suppose 
$$
L \in \mathcal{L}_{G_\sigma} \Leftrightarrow \Sigma := LDL^T \in P_{G_\sigma}. 
$$

\noindent
Then $G$ is a decomposable graph and $\sigma$ corresponds to a perfect vertex elimination 
scheme for $G$. 
\end{lemma}

\noindent
Hence, Lemma \ref{dcmplsprpl} and Lemma \ref{cnvrsdmpgr} characterize a decomposable graph 
$G$ and a perfect vertex elimination scheme $\sigma$ for $G$ in terms of the preservation of 
zeros in the modified Cholesky decomposition of matrices in $P_{G_\sigma}$. These 
characterizations of decomposable graphs and orderings of vertices of $G$ has proven to be 
tremendously useful for working with sparse positive definite matrices in probability and 
statistics (see \cite{khrrjtmwdg, ltcmssmwdg, roveratocs,rajaratnam08, guillot12}). Another class of graphs that is 
also highly useful in this context is the class of co-chordal graphs or homogeneous graphs 
(see \cite{adsnwjrwhc, khrrjtmcvt, khrrjtmcwp, khrrjtmwdg, ltcmssmwdg}). Yet 
characterizations of homogeneous graphs, similar to the above for decomposable graphs, are 
not available. These characterizations are the subject of the rest of the paper.

\subsection{Homogeneous graphs} \label{hmgnsgphpy}

\noindent
A graph $G = (V,E)$ is defined to be co-chordal or homogeneous if for all $v,v'$ such that 
$(v,v') \in E$, either 
$$
\{u: \; u = v' \mbox{ or } (u,v') \in E\} \subseteq \{u: \; u = v \mbox{ or } (u,v) \in E\}, 
$$

\noindent
or 
$$
\{u: \; u = v \mbox{ or } (u,v) \in E\} \subseteq \{u: \; u = v' \mbox{ or } (u,v') \in E\}. 
$$

\begin{figure}
\begin{center}
\includegraphics[width=3in,height=2in]{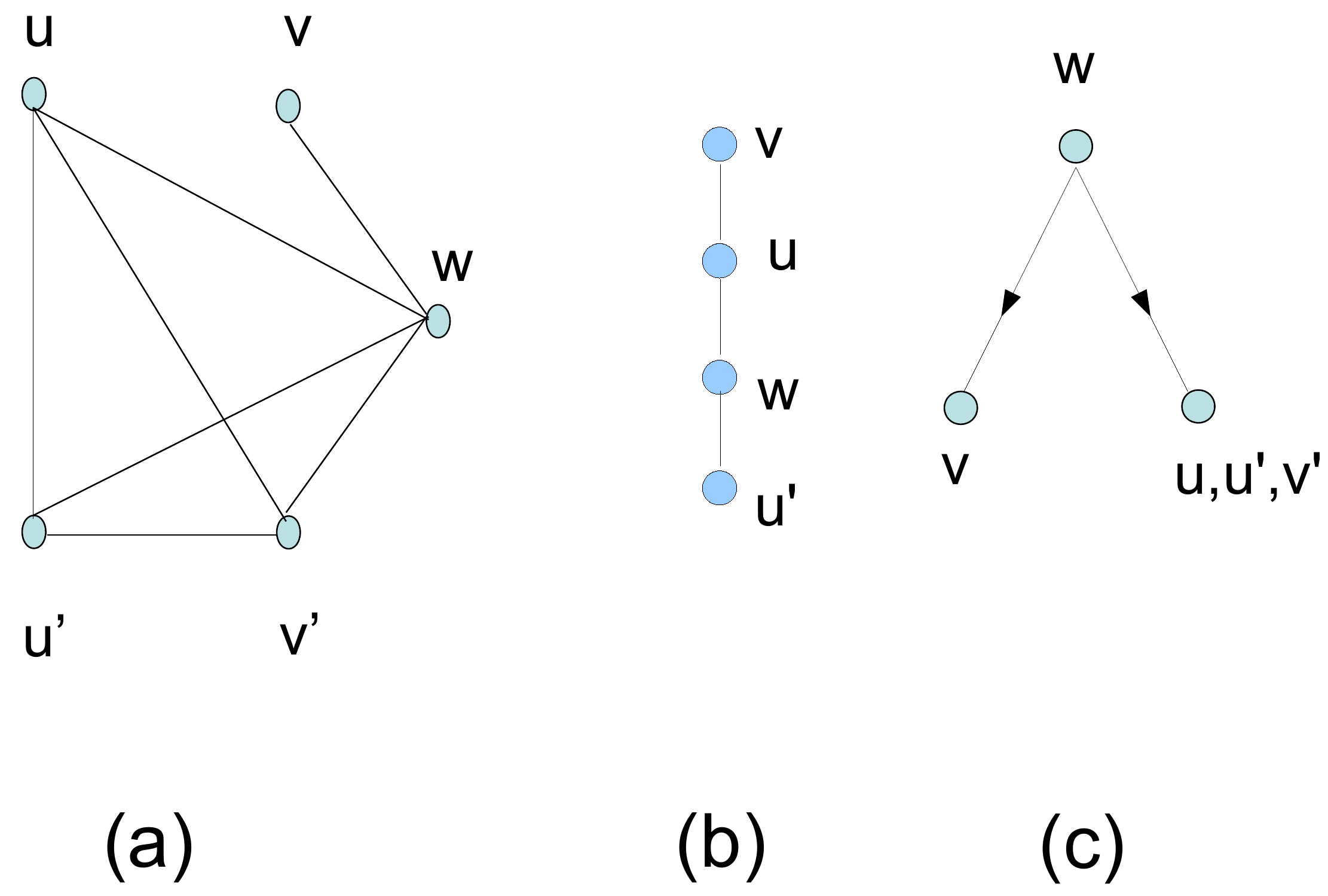}
\end{center}
\caption{(a) A homogeneous graph, (b) A non-homogeneous graph which is decomposable, and 
(c) The Hasse tree corresponding to the homogeneous graph in (a).}
\label{figure2}
\end{figure}

\noindent
Equivalently, a graph $G$ is said to be homogeneous if it is decomposable and does not 
contain the graph $\stackrel{\mbox{\tiny{1}}}{\bullet} - \stackrel{\mbox{\tiny{2}}}{\bullet} 
- \stackrel{\mbox{\tiny{3}}}{\bullet} - \stackrel{\mbox{\tiny{4}}}{\bullet}$, denoted by 
$A_4$, as an induced subgraph. See Figure \ref{figure2} for an example of a homogeneous 
graph, and a non-homogeneous graph which is decomposable. Connected homogeneous graphs 
have an equivalent representation in terms of directed rooted trees, called 
{\it Hasse diagrams}. The reader is referred to \cite{ltcmssmwdg} for a detailed account 
of the properties of homogeneous graphs. We write $v \rightarrow w$ whenever 
$$
\{u: \; u = w \mbox{ or } (u,w) \in E\} \subseteq \{u: \; u = v \mbox{ or } (u,v) \in E\}. 
$$

\noindent
Now denote by $R$ the equivalence relation on $V$ defined by 
$$
u R v \Leftrightarrow u \rightarrow v \mbox{ and } v \rightarrow u. 
$$

\indent
Let $\bar{v}$ denote the equivalence class in $V/R$ containing $v$. The Hasse diagram 
of $G$ is defined as a directed graph with vertex set $V_H = V/R = \{\bar{v}: v \in V\}$ 
and edge set $E_H$ consisting of directed edges with $(\bar{u}, \bar{v}) \in E_H$ for 
$\bar{u} \neq \bar{v}$ if the following holds: $u \rightarrow v$ and $\nexists v'$ such 
that $u \rightarrow v' \rightarrow v, \bar{v'} \neq \bar{u}, \bar{v'} \neq \bar{v}$. 

\indent
If $G$ is a connected homogeneous graph, then the Hasse diagram described above is a 
directed rooted tree such that the number of children of a vertex is never equal to one. 
It was proved in \cite{ltcmssmwdg} that there is a one-to-one correspondence between the 
set of connected homogeneous graphs and the set of directed rooted trees with vertices 
weighted by positive integers ($w(\bar{u}) = |\bar{u}|$), such that no vertex has exactly 
one child. If $u \rightarrow v$ and $\bar{u} \neq \bar{v}$, we say that $u$ is an ancestor 
of $v$ in the Hasse tree of $G$. It is easily seen that if $G$ is a disconnected 
homogeneous graph, then  each connected component of $G$ gives rise to a Hasse tree. If 
$\bar{u} = \bar{v}$, we say that $u$ is a twin of $v$ in the Hasse tree of $G$. 

\indent
A subclass of orderings associated with a homogeneous graph, which will be used in 
subsequent analysis, is defined as follows. 
\begin{definition}
If $G = (V,E)$ is a homogeneous graph, then an ordering $\sigma$ of $V$ is defined to be a 
Hasse tree based elimination scheme for $G$ if for every pair of vertices $u,v$, the 
following holds. 
$$
u \rightarrow v, \bar{u} \neq \bar{v} \Rightarrow \sigma(u) > \sigma(v). 
$$

\noindent
Alternatively, if $\bar{u}$ is an ancestor of $\bar{v}$ in the Hasse diagram of $G$, then 
$\sigma(u) > \sigma(v)$. 
\end{definition}

\noindent
The lemma below follows easily from the definition of homogeneous graphs. 
\begin{lemma} \label{hmgnsunion}
(a) If $G_i = (V_i, E_i)$ is a homogeneous graph for every $1 \leq i \leq n$, and $V_i$ and 
$V_j$ are disjoint for every $1 \leq i \neq j \leq n$, then $G = \left( \cup_{i=1}^n V_i, 
\cup_{i=1}^n E_i \right)$ is also a homogeneous graph. Conversely, if $G = (V,E)$ is a 
homogeneous graph, then any disjoint connected component of $G$ is also a homogeneous 
graph. 

\noindent
(b) If $G = (V,E)$ is a connected homogeneous graph, $|V| = m$, and $\sigma$ is a Hasse 
tree based elimination scheme for $G$, then the equivalence class of $\sigma^{-1} (m)$ lies 
at the root of the Hasse tree of $G$. 
\end{lemma}

\newtheorem{example}{Example}
\begin{example}
Consider the homogeneous graph $G$ in Figure \ref{figure2} {\it (a)} and the corresponding 
Hasse tree in Figure \ref{figure2} {\it (c)}. A Hasse tree based elimination scheme $\sigma$ 
for the homogeneous graph $G$ is given by $\sigma(w) = 5, \sigma(v) = 4, \sigma(v') = 3, 
\sigma(u') = 2, \sigma(u) = 1$. Note that a homogeneous graph is also a decomposable graph, 
and a Hasse tree based elimination scheme is also a perfect vertex elimination scheme. 
However, every perfect vertex elimination scheme for a homogeneous graph may not necessarily 
be a Hasse tree based elimination scheme. For the homogeneous graph $G$ in 
Figure \ref{figure2} {\it (a)}, the ordering $\sigma$ given by $\sigma(v') = 5, \sigma(w) = 
4, \sigma(u') = 3, \sigma(u) = 2, \sigma(v) = 1$ is a perfect vertex elimination scheme, but 
not a Hasse tree based elimination scheme, since $w \rightarrow v', \bar{w} \neq \bar{v'}$ 
but $\sigma(w) = 4 < \sigma(v') = 5$. 
\end{example}

\section{Characterization in terms of sparse matrix decompositions} \label{first_characterization}

\noindent
We now provide the first characterization of homogeneous graphs that yields a parallel 
result to that of Paulsen et al. \cite{plsnpwrsth} for decomposable graphs. We note that antecedents of the results in Paulsen et al. \cite{plsnpwrsth} were given in \cite{dym81, grone84, agler88}.

\begin{lemma}[Khare and Rajaratnam \cite{khrrjtmwdg}] \label{hmgnshtbes}
Let $G = (V,E)$ be a homogeneous graph, and $\sigma$ an ordering of $V$ which corresponds 
to a Hasse tree based elimination scheme for $G$. Then for any positive definite matrix 
$\Sigma$ with modified Cholesky decomposition given by $\Sigma = LDL^T$, the following 
holds. 
$$
\Sigma \in P_{G_\sigma} \Leftrightarrow L \in L_{G_\sigma} \Leftrightarrow L^{-1} \in 
L_{G_\sigma}. 
$$
\end{lemma}

\noindent
A detailed constructive proof is given in \cite{khrrjtmcwp}. A proof in a more general 
context can also be found in \cite{adsnwjrwhc} and \cite{trgpapcnhr}. One of the main 
results of this paper is the converse of Lemma \ref{hmgnshtbes}. 
\newtheorem{prop}{Proposition}
\begin{prop} \label{cnvrshgsht}
Let $G = (V,E)$ be a graph, $\sigma$ be an ordering of $V$, and $D$ be an arbitrary diagonal 
matrix with positive diagonal entries. Suppose 
$$
L \in \mathcal{L}_{G_\sigma} \Leftrightarrow L^{-1} \in \mathcal{L}_{G_\sigma} 
\Leftrightarrow \Sigma := LDL^T \in P_{G_\sigma}. 
$$

\noindent
Then $G$ is a homogeneous graph and $\sigma$ corresponds to a Hasse tree based elimination 
scheme for $G$. 
\end{prop}

\noindent
{\it Proof:} We proceed by induction and prove the result in a series of claims. 
\vspace{0.3cm}

\noindent
{\bf Claim 1:} The result holds for $|V| = 3$. 

\noindent
{\it Proof of Claim 1}: Let $V = \{u,v,w\}$. We consider two cases. 

\begin{figure}
\begin{center}
\includegraphics[width=4in,height=1in]{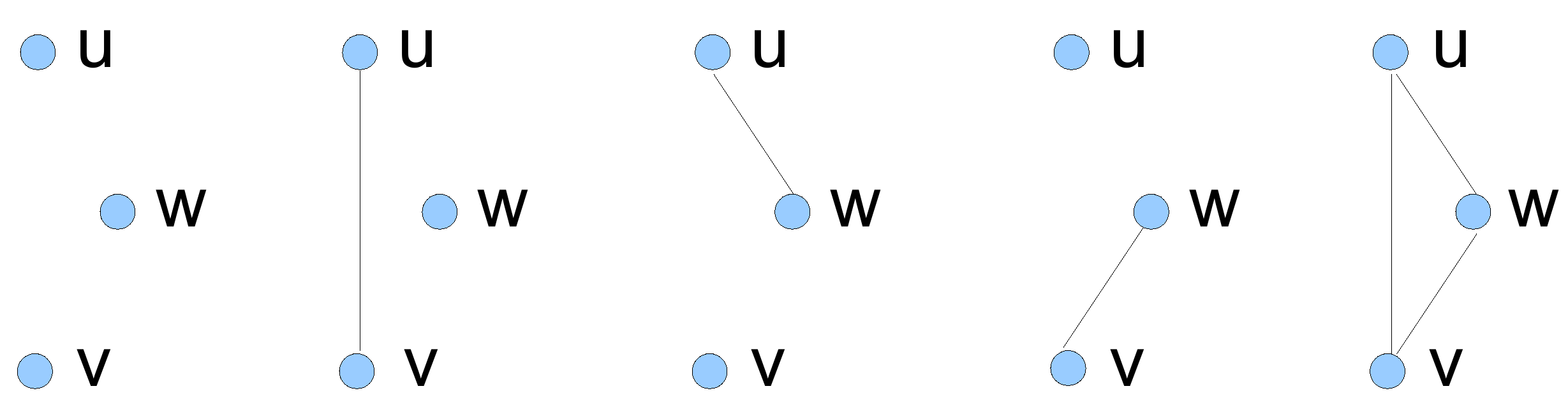}
\end{center}
\caption{Case I with $|V| = 3$ for Proposition \ref{cnvrshgsht}}
\label{figure3}
\end{figure}

\vspace{0.15cm}
\noindent
{\bf Case I:} $E = \phi, \{(u,v)\}, \{(u,w)\}, \{(v,w)\} \mbox{ or } \{(u,v),(u,w),(v,w)\}$. See Figure \ref{figure3}. 

\noindent
$G$ is a homogeneous graph in every case. Also, each disjoint connected component is a 
complete graph, which means that every ordering corresponds to a Hasse tree based 
elimination scheme. Hence, the result holds vacuously. 

\begin{figure}
\begin{center}
\includegraphics[width=4in,height=1in]{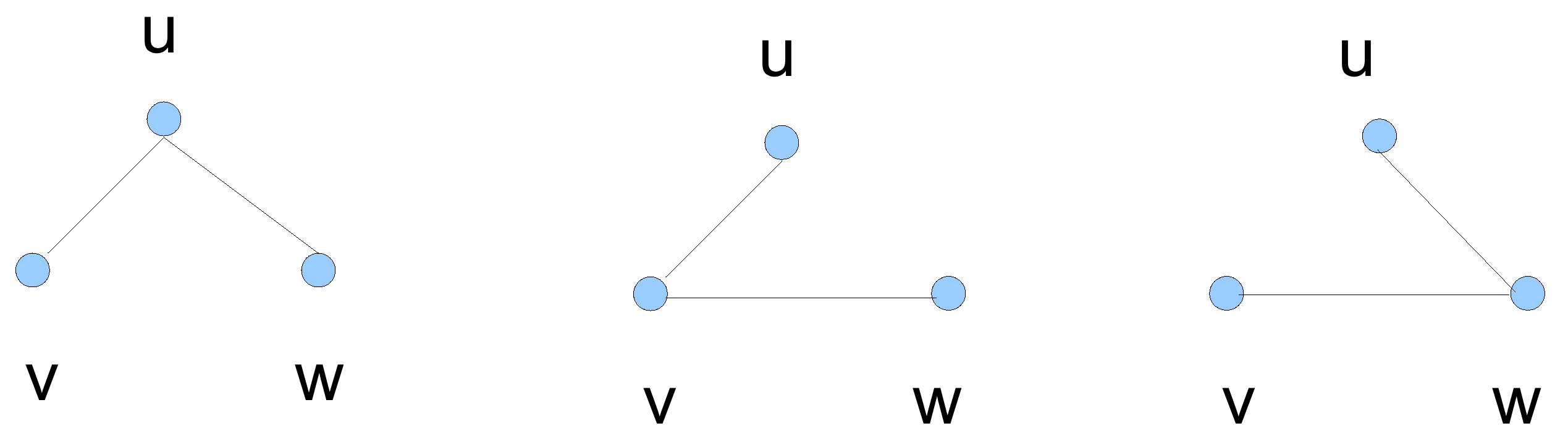}
\end{center}
\caption{Case II with $|V| = 3$ for Proposition \ref{cnvrshgsht}}
\label{figure4}
\end{figure}

\vspace{0.15cm}
\noindent
{\bf Case II:} $E = \{(u,v),(v,w)\}, \{(u,w),(v,w)\} \mbox{ or } \{(u,v),(u,w)\}$. See Figure \ref{figure4}.

\noindent
Let us first consider the case $E = \{(u,v),(v,w)\}$. Note that $G$ is a homogeneous graph. 
It remains to be shown that $\sigma$ is a Hasse tree based elimination scheme. Now if 
$\sigma(v) = 1$, and 
$$
L = \left( \begin{matrix}
1 & 0 & 0 \cr
1 & 1 & 0 \cr
1 & 0 & 1 
\end{matrix}
\right) \in \mathcal{L}_{G_\sigma}, 
$$

\noindent
then $\Sigma_{32} = (LDL^T)_{32} = d_{11} \neq 0$. Hence, $\Sigma \notin P_{G_\sigma}$, 
yielding a contradiction. Similarly, if $\sigma(v) = 2$, and 
$$
L = \left( \begin{matrix}
1 & 0 & 0 \cr
1 & 1 & 0 \cr
0 & 1 & 1 
\end{matrix}
\right) \in \mathcal{L}_{G_\sigma}, 
$$

\noindent
then $L^{-1}_{31} = 1 \neq 0$. Hence, $L^{-1} \notin \mathcal{L}_{G_\sigma}$, once more 
yielding a contradiction to the assumptions in the proposition. Hence $\sigma(v) = 3$. Note 
that $v \rightarrow u, v \rightarrow w$ and $\bar{v} \neq \bar{u}, \bar{v} \neq \bar{w}$. 
Hence, any ordering $\sigma$ such that $\sigma(v) = 3$ is a Hasse tree based elimination 
scheme. The other cases when $E = \{(u,w), (v,w)\}$ and $E = \{(u,v), (u,w)\}$ follow by 
symmetry. Hence, the result for $|V| = 3$ holds true. 

\vspace{0.1in}

\indent
As mentioned earlier, we shall use an induction argument on the number of vertices to prove 
the result. Suppose now that the result holds true for all graphs with $m-1$ vertices. Let 
$G = (V,E)$ be a graph with $|V| = m$, and $\sigma$ be an ordering of $V$ for which 
$$
L \in \mathcal{L}_{G_\sigma} \Leftrightarrow L^{-1} \in L_{G_\sigma} \Leftrightarrow 
\Sigma := LDL^T \in P_{G_\sigma}, 
$$

\noindent
for an arbitrary diagonal matrix $D$ with positive diagonal entries. We need to show two 
results: (i) $G$ is homogeneous and (ii) the ordering $\sigma$ is a Hasse tree based 
elimination scheme. 

\indent
Let $G'$ be the subgraph induced by $G$ on the set of vertices $V \setminus \{\sigma^{-1} 
(m)\}$, and let $\sigma'$ be the restriction of $\sigma$ on $V \setminus \{\sigma^{-1} 
(m)\}$. Note that $G'$ together with the ordering $\sigma'$ is none other than $G$ with the 
ordering $\sigma$ (or $G_\sigma$), but with the highest labeled vertex removed. \\

\noindent {\bf Claim 2}: 
$$
L^* \in \mathcal{L}_{G'_{\sigma'}} \Leftrightarrow (L^*)^{-1} \in \mathcal{L}_{G'_{\sigma'}} 
\Leftrightarrow \Sigma^* = L^*D^*(L^*)^T \in P_{G'_{\sigma'}}. 
$$

where $D^*$ is the upper $(m-1) \times(m-1)$ principal submatrix of $D$. 

\vspace{0.1in}

\noindent
{\it Proof of Claim 2}: Let $L^* \in \mathcal{L}_{G'_{\sigma'}}$. Then 
\begin{eqnarray*}
& & L := \left( \begin{matrix}
L^* & {\bf 0} \cr
{\bf 0}^T & 1 
\end{matrix}
\right) \in \mathcal{L}_{G_\sigma}\\
&\Rightarrow& \left( \begin{matrix}
L^* & {\bf 0} \cr
{\bf 0}^T & 1 
\end{matrix}
\right)^{-1} = \left( \begin{matrix}
(L^*)^{-1} & {\bf 0} \cr
{\bf 0}^T & 1 
\end{matrix}
\right) \in \mathcal{L}_{G_\sigma}\\
&\Rightarrow& (L^*)^{-1} \in \mathcal{L}_{G'_{\sigma'}}. 
\end{eqnarray*}

\noindent
By a similar argument $(L^*)^{-1} \in \mathcal{L}_{G'_{\sigma'}} \Rightarrow L^* \in 
\mathcal{L}_{G'_{\sigma'}}$. Hence $(L^*)^{-1} \in \mathcal{L}_{G'_{\sigma'}} 
\Leftrightarrow L^* \in \mathcal{L}_{G'_{\sigma'}}$. 

\noindent Note that, 
\begin{eqnarray*}
& & L^* \in \mathcal{L}_{G'_{\sigma'}}\\
&\Leftrightarrow& L = \left( \begin{matrix}
L^* & {\bf 0} \cr
{\bf 0}^T & 1 
\end{matrix}
\right) \in \mathcal{L}_{G_\sigma}\\
&\Leftrightarrow& \Sigma = LDL^T = \left( \begin{matrix}
L^*D^*(L^*)^T & {\bf 0} \cr
{\bf 0}^T & D_{mm} 
\end{matrix}
\right) \in P_{G_\sigma}\\
&\Leftrightarrow& \Sigma^* := L^*D^*(L^*)^T \in P_{G'_{\sigma'}}. 
\end{eqnarray*}

\noindent
Hence, we have now established that 
$$
L^* \in \mathcal{L}_{G'_{\sigma'}} \Leftrightarrow (L^*)^{-1} \in \mathcal{L}_{G'_{\sigma'}} 
\Leftrightarrow \Sigma^* = L^*D^*(L^*)^T \in P_{G'_{\sigma'}}. 
$$

\vspace{0.1in}

\noindent
By the induction hypothesis, it follows that $G'$ is a homogeneous graph and $\sigma'$ 
corresponds to a Hasse tree based elimination scheme for $G'$, i.e., 
\begin{equation} \label{indtnhasod}
\sigma'(v) = \sigma(v) < \sigma'(u) = \sigma(u) \mbox{ when } u \rightarrow v, \bar{u} 
\neq \bar{v}, \; \forall u,v \in V \setminus \{\sigma^{-1}(m)\}. 
\end{equation}

\noindent
{\bf Claim 3}: $G$ is a homogeneous graph and $\sigma$ is a Hasse tree based elimination 
scheme. 

\vspace{0.1in}

\noindent
{\it Proof of Claim 3:}  Now let $V' = \cup_{i=1}^k V_i$, where $V_i$ is the vertex set 
corresponding to the $i^{th}$ disjoint connected component of $G'$. 

\indent
Suppose $(\sigma^{-1} (m), u) \notin E$ for each $u \in V \setminus \{\sigma^{-1}(m)\}$, 
i.e., the vertex $\sigma^{-1}(m)$ is disconnected from the graph $G'$. Then by 
Lemma \ref{hmgnsunion}, the graph $G$ is a homogeneous graph with $V = \left( \cup_{i=1}^k 
V_i \right) \cup \{\sigma^{-1}(m)\}$ being the disjoint partition of the vertices 
corresponding to its disjoint connected components. Also, from (\ref{indtnhasod}) and the 
fact that $\sigma^{-1}(m)$ is disconnected from every vertex in $V \setminus 
\{\sigma^{-1}(m)\}$, it follows that $\sigma$ is a Hasse tree based elimination scheme for 
$G$. 

\indent
Suppose $(\sigma^{-1}(m), u) \in E$ for some $u \in V_i$. Let $v_i^* \in V_i$ be such that 
$\sigma (v_i^*) = \max_{v_i \in V_i} \sigma (v_i)$. Since $G'$ is a homogeneous graph, 
$\sigma$ restricted to $V \setminus \{\sigma^{-1}(m)\}$ is a Hasse tree based elimination 
scheme, and $V_i$ is the vertex set corresponding to a connected component of $G'$, it 
follows from Lemma \ref{hmgnsunion} that the equivalence class of $v_i^*$ lies at the top 
of the Hasse tree of $V_i$ in $G'$. We therefore deduce that $(v_i^*, v_i) \in E, \; 
\forall v_i \in V_i$. 

\indent
We proceed by claiming that $(\sigma^{-1}(m), v_i^*) \in E$. If $v_i^* = u$, it follows 
immediately. If $v_i^* \neq u$, then $m > \sigma(v_i^*) > \sigma(u)$. Suppose $L$ is 
defined by 
$$
L_{ij} = \begin{cases}
1 & \mbox{if } i = m, j = \sigma(u) \mbox{ or } i = \sigma(v_i^*), j = \sigma(u) 
\mbox{ or } i = j, \cr 
0 & \mbox{otherwise}. 
\end{cases}
$$

\noindent
Note that $L \in \mathcal{L}_{G_\sigma}$. If $\Sigma := LL^T$, then by assumption $\Sigma 
\in P_{G_\sigma}$, and 
$$
\Sigma_{m \sigma(v_i^*)} = L_{m \sigma(u)} L_{\sigma(v_i^*) \sigma(u)} + \sum_{v \in V_i, 
v \neq u} L_{m \sigma(v)} L_{\sigma(v_i^*) \sigma(v)} = 1. 
$$

\noindent
Hence, it follows that $(\sigma^{-1}(m), v_i^*) \in E$. Now let $v_i \in V_i, v_i \neq 
v_i^*$. We also now claim that $(\sigma^{-1}(m), v_i) \in E$. Note that $(v_i^*, v_i) \in 
E$ from the discussion above. Suppose $L$ is defined by 
$$
L_{ij} = \begin{cases}
1 & \mbox{if } i = m, j = \sigma(v_i^*) \mbox{ or } i = \sigma(v_i^*), j = \sigma(v_i) 
\mbox{ or } i = j, \cr
0 & \mbox{otherwise}. 
\end{cases}
$$

\noindent
First note that $L \in \mathcal{L}_{G_\sigma}$, and hence by assumption $L^{-1} \in 
\mathcal{L}_{G_\sigma}$. Since $L^{-1}_{m \sigma(v_i)} = 1$ (by using the inversion formula 
in Lemma \ref{inverselrt}), it follows that $(\sigma^{-1}(m), v_i) \in E$. Hence, we have 
established that if $(\sigma^{-1} (m), u) \in E$ for some $u \in V_i$, then $(\sigma^{-1} 
(m), v_i) \in E$ for every $v_i \in V_i$. 

\indent
Now let $V_{i_1}, V_{i_2}, \cdots, V_{i_p}$ be the components of $G'$ which share at least 
one edge with $\sigma^{-1}(m)$. Since the graph induced by $V_{i_r}$ on $G'$ is a connected 
homogeneous graph for every $1 \leq r \leq p$, and $\sigma^{-1} (m)$ is connected to every 
vertex in $V_{i_1}, V_{i_2}, \cdots, V_{i_p}$ by the argument above, the introduction of 
$\sigma^{-1} (m)$ does not give rise to any new $4$-cycle or $4$-path, due of the 
following reasoning: Consider an arbitrary collection of $4$ vertices in $V$. If all of them lie in $V_{i_r}$ 
for some $r$, and if $\sigma^{-1} (m)$ is not one of the vertices, then these 4 vertices cannot form a $4$-cycle or a $4$-path as the subgraph induced by $V_{i_r}$ on $G$ is a homogeneous graph. If none of the vertices is $\sigma^{-1} (m)$, and all of them do not lie in $V_{i_r}$ for some $r$, then the graph induced by these vertices on $G$ 
is a disconnected graph, which implies that the induced sub-graph cannot be a $4$-cycle or a $4$-path. Finally, if $\sigma^{-1} (m)$ is one of the vertices, and since it is connected to all the other three vertices, they cannot form an induced $4$-cycle or an induced $4$-path.   

\indent
It follows that the graph induced by $\{\sigma^{-1}(m)\} \cup \left( \cup_{r=1}^p V_{i_r} 
\right)$ on $G$ is a connected homogeneous graph. Moreover, since $\sigma^{-1} (m)$ is 
connected to every vertex in $V_{i_1}, V_{i_2}, \cdots, V_{i_p}$, its equivalence class has 
to lie at the root of the corresponding Hasse tree. Note that the disjoint connected 
components of $G'$ other than $V_{i_1}, V_{i_2}, \cdots, V_{i_p}$ are also connected 
homogeneous graphs. It follows that $G$ is a homogeneous graph with disjoint connected 
components $\{\sigma^{-1}(m)\} \cup \left( \cup_{r=1}^p V_{i_r} \right)$ and $V_t, \; t \neq 
i_1, i_2, \cdots, i_p$. Note that $\sigma'$ (which is the restriction of $\sigma$ to $G'$) corresponds to a Hasse tree based elimination scheme for $G'$, and that $\sigma(u) < m$ whenever $u \neq \sigma^{-1} (m)$. Hence, 
$\sigma(u) < m$ whenever $\sigma^{-1} (m) \rightarrow u, \overline{\sigma^{-1} (m)} \neq 
\bar{u}$. Also, since $\sigma^{-1} (m)$ is at the top of the Hasse tree in its connected 
component, there does not exist $u \in V \setminus \{\sigma^{-1} (m)\}$ such that $u 
\rightarrow \sigma^{-1} (m)$. This leads us to conclude that $\sigma$ is a Hasse tree 
based elimination scheme for $G$. Hence the result is proved. \hfill$\Box$

\vspace{0.1in}

{\it Remark:} A useful alternative probabilistic characterization of homogeneous graphs can be found
in \cite{pearlwrmth, drtonrchsn}. This probabilistic result essentially states that $G$ is homogeneous iff ``$G$ is Markov equivalent to a directed acyclic graph(DAG)". In contrast, the characterization proved in this section is
algebraic in nature, and is therefore different from the probabilistic characterization. The algebraic characterization above can be established directly starting from the probabilistic characterization mentioned above,
by using the notion of ``d-separation". The proof however is non-trivial and does not seem to offer a simplification over the first principles proof provided here.

\vspace{0.1in}

\noindent
We now give a series of examples to illustrate the necessity of the assumptions in the 
characterization discussed above. 
\begin{example}
Consider the homogeneous graph $G$ in Figure \ref{figure2} {\it (a)}. Let 
$\sigma$ be a Hasse tree based elimination scheme defined by $\sigma(w) = 5, \sigma(v) = 4, 
\sigma(v') = 3, \sigma(u') = 2, \sigma(u) = 1$. Let 
$$
L = \left( \begin{matrix}
1 & 0 & 0 & 0 & 0 \cr
1 & 1 & 0 & 0 & 0 \cr
1 & 1 & 1 & 0 & 0 \cr
0 & 0 & 0 & 1 & 0 \cr
1 & 1 & 1 & 1 & 1 
\end{matrix} \right) \in \mathcal{L}_{G_\sigma}. 
$$

\noindent
Then, 
$$
L^{-1} = \left( \begin{matrix}
1 & 0 & 0 & 0 & 0 \cr
-1 & 1 & 0 & 0 & 0 \cr
0 & -1 & 1 & 0 & 0 \cr
0 & 0 & 0 & 1 & 0 \cr
0 & 0 & -1 & -1 & 1 
\end{matrix} \right) \in \mathcal{L}_{G_\sigma}, \mbox{ and } 
\Sigma = LL^T = \left( \begin{matrix}
1 & 1 & 1 & 0 & 1 \cr
1 & 2 & 2 & 0 & 2 \cr
1 & 2 & 3 & 0 & 3 \cr
0 & 0 & 0 & 1 & 1 \cr
1 & 2 & 3 & 1 & 5 
\end{matrix} \right) \in P_{G_\sigma}. 
$$

\noindent
Now consider $\sigma$ which is a perfect vertex elimination scheme, but not a Hasse tree 
based elimination scheme, given by $\sigma(v') = 5, \sigma(w) = 4, \sigma(u') = 3, \sigma(u) 
= 2, \sigma(v) = 1$. Then 
$$
L = \left( \begin{matrix}
1 & 0 & 0 & 0 & 0 \cr
0 & 1 & 0 & 0 & 0 \cr
0 & 1 & 1 & 0 & 0 \cr
1 & 1 & 1 & 1 & 0 \cr
0 & 1 & 1 & 1 & 1 
\end{matrix} \right) \in \mathcal{L}_{G_\sigma}, \mbox{ but } 
L^{-1} = \left( \begin{matrix}
1 & 0 & 0 & 0 & 0 \cr
0 & 1 & 0 & 0 & 0 \cr
0 & -1 & 1 & 0 & 0 \cr
-1 & 0 & -1 & 1 & 0 \cr
1 & 0 & 0 & -1 & 1 
\end{matrix} \right) \notin \mathcal{L}_{G_\sigma}. 
$$

\noindent
It can be verified that $\Sigma = LL^T \in P_{G_\sigma}$. Now let $\sigma$ be given by 
$\sigma(v) = 5, \sigma(u') = 4, \sigma(v') = 3, \sigma(w) = 2, \sigma(u) = 1$. Then, 
$\sigma$ is not a perfect vertex elimination scheme, and 
$$
L = \left( \begin{matrix}
1 & 0 & 0 & 0 & 0 \cr
1 & 1 & 0 & 0 & 0 \cr
1 & 1 & 1 & 0 & 0 \cr
1 & 1 & 1 & 1 & 0 \cr
0 & 1 & 0 & 0 & 1 
\end{matrix} \right) \in \mathcal{L}_{G_\sigma}, \mbox{ but } 
\Sigma = LL^{T} = \left( \begin{matrix}
1 & 1 & 1 & 1 & 0 \cr
1 & 2 & 2 & 2 & 1 \cr
1 & 2 & 3 & 3 & 1 \cr
1 & 2 & 3 & 4 & 1 \cr
0 & 1 & 1 & 1 & 2 
\end{matrix} \right) \notin P_{G_\sigma}. 
$$

\noindent
Now consider the non-homogeneous graph $G$ in Figure \ref{figure2} {\it (b)}. Note that $G$ 
is however a decomposable graph. The ordering $\sigma$ given by $\sigma(u') = 4, \sigma(w) 
= 3, \sigma(u) = 2, \sigma(v) = 1$ is a perfect vertex elimination scheme. However, 
$$
L = \left( \begin{matrix}
1 & 0 & 0 & 0 \cr
1 & 1 & 0 & 0 \cr
0 & 1 & 1 & 0 \cr
0 & 0 & 1 & 1 
\end{matrix} \right) \in \mathcal{L}_{G_\sigma}, \mbox{ but } 
L^{-1} = \left( \begin{matrix}
1 & 0 & 0 & 0 \cr
-1 & 1 & 0 & 0 \cr
1 & -1 & 1 & 0 \cr
-1 & 1 & -1 & 1 
\end{matrix} \right) \notin \mathcal{L}_{G_\sigma}. 
$$
\end{example}

\section{Characterization in terms of determinants} \label{second_characterization}

\noindent
We now give a second characterization of homogeneous graphs with vertex orderings 
corresponding to Hasse tree based elimination schemes. Let us first establish some 
notation, that shall be used throughout this section. If $A \in \mathbb{M}_n$ and $M, M^* 
\subseteq \{1,2 \cdots, n\}$, then 
$$
A_M := ((A_{ij}))_{i,j \in M}, \; \; A_{MM^*} := ((A_{ij}))_{i \in M, j \in M^*}. 
$$

\noindent
The proposition below and its converse, stated and proved subsequently, provide the second 
characterization of homogeneous graphs. 
\begin{prop} \label{hmgnsdtmnt}
Let $G = (V,E)$ be a homogeneous graph, and $\sigma$ an ordering of $V$ which corresponds to 
a Hasse tree based elimination scheme for $G$. Let $\Sigma \in P_{G_\sigma}$, and $\Sigma = 
LDL^T$ denote its Cholesky decomposition. Then, for any maximal clique $C$, 
$$
\left| (\Sigma^{-1})_{\sigma(C)} \right| = \prod_{i \in \sigma(C)} \frac{1}{D_{ii}}. 
$$
\end{prop}

\noindent {\it Proof:}  Let $C \subseteq V$ be a maximal clique in $G$, where $C = \{u_1, u_2, \cdots, 
u_r\}$, with $\sigma(u_1) > \sigma(u_2) > \cdots > \sigma(u_r)$. First note that 
\begin{equation} \label{mxmlclqdet}
(\Sigma^{-1})_{\sigma(C)} = \left[ (L^{-1})_{\sigma(V) \sigma(C)} \right]^T D^{-1} \left[ 
(L^{-1})_{\sigma(V) \sigma(C)} \right]. 
\end{equation}

\noindent
We will prove that the determinant of the RHS of (\ref{mxmlclqdet}) equals the determinant 
of \\ $\left[ (L^{-1})^T_{\sigma(C)} \right] D^{-1}_{\sigma(C)} \left[ (L^{-1})_{\sigma(C)} 
\right]$, and the result will follow. 

\indent
We start by first showing that $L^{-1}_{\sigma(w) \sigma(u_i)} = 0$ when $w \notin C$ for 
$i = 1,2, \cdots, r$. Note that $\sigma(u_i) > \sigma(w), \; L^{-1}_{\sigma(w) \sigma(u_i)} 
= 0$, as $L^{-1}$ is a lower triangular matrix. Now let $\sigma(u_i) < \sigma(w)$. Suppose to 
the contrary that $L^{-1}_{\sigma(w) \sigma(u_i)} \neq 0$. Since $L^{-1} \in 
\mathcal{L}_{G_\sigma}$ by Lemma \ref{hmgnshtbes}, we get $(w, u_i) \in E$. Hence, $w$ is 
an ancestor or twin of $u_i$ in the Hasse tree of $G$. Now by the very definition of a 
homogeneous graph, every vertex sharing an edge with $u_i$ also shares an edge with $w$. 
Hence, $(w,u_j) \in E$ for $j = 1,2, \cdots, r$, which gives a contradiction to the 
maximality of $C$. Hence we conclude that $L^{-1}_{\sigma(w) \sigma(u_i)} = 0$ when $w 
\notin C$ for $i = 1,2, \cdots, r$. 

\indent
Now using the Cauchy-Binet identity in (\ref{mxmlclqdet}), 
\begin{eqnarray*}
\left| (\Sigma^{-1})_{\sigma(C)} \right| 
&=& \left| \left[ (L^{-1})_{\sigma(V) \sigma(C)} \right]^T D^{-1} \left[ (L^{-1})_{\sigma(V) 
\sigma(C)} \right] \right|\\
&=& \sum_{A \subseteq V, |A| = r} \left| \left[ (L^{-1})_{\sigma(A) \sigma(C)} \right]^T 
D^{-1}_{\sigma(A)} \left[ (L^{-1})_{\sigma(A) \sigma(C)} \right] \right|. 
\end{eqnarray*}

\noindent
Note that if $A \subseteq V, |A| = r$, and $A \neq C$, then there exists $w$ such that $w 
\in A$ but $w \notin C$. Hence, from the argument above, $L^{-1}_{\sigma(w) \sigma(u_i)} = 
0$ for $i = 1,2, \cdots, r$, and for such $A \neq C$, 
$$
\left| \left[ (L^{-1})_{\sigma(A) \sigma(C)} \right]^T D^{-1}_{\sigma(A)} \left[ 
(L^{-1})_{\sigma(A) \sigma(C)} \right] \right| = \left| \left[ (L^{-1})_{\sigma(A) 
\sigma(C)} \right]^T \right| \left| D^{-1}_{\sigma(A)} \right| \left| \left[ 
(L^{-1})_{\sigma(A) \sigma(C)} \right] \right| = 0, 
$$

\noindent
since one row in the matrix $(L^{-1})_{\sigma(A) \sigma(C)}$ is zero. Therefore the only 
non-zero summand in the Cauchy-Binet formula is when $A = C$. Hence 
$$
\left| (\Sigma^{-1})_{\sigma(C)} \right| = \left| \left[ (L^{-1})_{\sigma(C)} \right]^T 
\right| \left| D^{-1}_{\sigma(C)} \right| \left| \left[ (L^{-1})_{\sigma(C)} \right] \right| 
= \prod_{i \in \sigma(C)} \frac{1}{D_{ii}}, 
$$

\noindent
where the last equality follows from the fact that $(L^{-1})_{\sigma(C)}$ is a lower 
triangular matrix with all diagonal entries equal to one (and therefore has determinant 
one) , and $D^{-1}_{\sigma(C)}$ is a diagonal matrix. Hence the result is proved. 
\hfill$\Box$ 

\vspace{0.1in}

\noindent
We now proceed to prove the following lemma required in the proof of the converse of 
Proposition \ref{hmgnsdtmnt}. 
\begin{lemma} \label{cyclptharg}
Let $G = (V,E)$ be a $4$-cycle or $4$-path, and let $\sigma$ be an ordering of $V$. Then, 
irrespective of the way $\sigma$ orders the vertices of the $4$-cycle or the $4$-path, 
there exist $u,v,w \in V$ such that $(u,v), (v,w) \in E, \; (u,w) \notin E$, and $\sigma(v) 
< \sigma(u) < \sigma(w)$ or $\sigma(u) < \sigma(v) < \sigma(w)$. 
\end{lemma}

\noindent
{\it Proof:}  (i) Let $G$ be a $4$-cycle. Recall that $u,v \in V$ are said to be 
{\it neighbors} in $G$ if $(u,v) \in E$. Consider the two neighbors of $v := \sigma^{-1} 
(1)$. Let $u$ denote the neighbor with the smaller $\sigma$-value, and $w$ denote the 
remaining neighbor. Note that $(u,v), (v,w) \in E$, but $(u,w) \notin E$. Also, $\sigma(v) 
= 1 < \sigma(u) < \sigma (w)$. 

\noindent
(ii) Let $G$ be a $4$-path. We consider three possibilities which are exhaustive, and in 
each case show the existence of three vertices with the required properties. 
\begin{enumerate}
\item[Case I] {\bf $\sigma^{-1} (1)$ has two neighbors}: Let $v := \sigma^{-1} (1)$. In this case, let $u$ denote the 
neighbor with the smaller $\sigma$-value, and $w$ denote the remaining neighbor. 
Hence, $\sigma (v) = 1 < \sigma(u) < \sigma(w)$. 
\item[Case II]{\bf $\sigma^{-1} (1)$ has one neighbor, and $\sigma^{-1} (2)$ has two neighbors}: 
Let $v := \sigma^{-1} (2)$. If one of the two neighbors of $v = \sigma^{-1} (2)$ is $u = \sigma^{-1} (1)$, denote the 
remaining neighbor by $w$, and observe that $\sigma(w)$ is equal to $3$ or $4$. Hence, 
$\sigma(u) = 1 < \sigma(v) = 2 < \sigma(w)$. If the neighbors of $v = \sigma^{-1} (2)$ are 
$u = \sigma^{-1} (3)$ and $w = \sigma^{-1} (4)$, then $\sigma(v) = 2 < \sigma(u) = 3 < 
\sigma(w) = 4$. 
\item[Case III] {\bf $\sigma^{-1} (1)$ and $\sigma^{-1} (2)$ both have one neighbor}: In this case, 
$v := \sigma^{-1} (3)$ has two neighbors, one of which has to be $w = \sigma^{-1} (4)$. 
Let $u$ be the remaining neighbor and observe that $\sigma(u)$ is equal to $1$ or $2$. 
Hence, $\sigma(u) < \sigma(v) = 3 < \sigma(w) = 4$. 
\end{enumerate}

\vspace{0.1in}

\noindent
We now establish the converse of Proposition \ref{hmgnsdtmnt}. 
\begin{prop} \label{cnvrshgsdt}
Let $G = (V,E)$ be a graph, and $\sigma$ be an ordering of $V$. Now if $G$ is not a 
homogeneous graph, or if $G$ is a homogeneous graph and $\sigma$ does not correspond to a 
Hasse tree based elimination scheme for $G$, then there exists a maximal clique $C$, and 
$\Sigma \in P_{G_\sigma}$ such that 
$$
\left| (\Sigma^{-1})_{\sigma(C)} \right| \neq \prod_{i \in \sigma(C)} \frac{1}{D_{ii}}, 
$$

\noindent
where $\Sigma = LDL^T$ denotes the modified Cholesky decomposition of $\Sigma$. 
\end{prop}

\noindent {\it Proof of Proposition \ref{cnvrshgsdt}:} We shall prove the result for each of the two possible cases. 
\vspace{0.3cm}

\noindent
{\bf Case I}: $G$ is not a homogeneous graph. 

\noindent
As the graph $G$ is not homogneous, it contains a $4$-cycle or a $4$-path. If $G$ contains 
a $4$-cycle or a $4$-path, by Lemma \ref{cyclptharg}, there exist $u,v,w \in V$ such that $(u,v), (v,w) 
\in E, \; (u,w) \notin E$, and $\sigma(v) < \sigma(u) < \sigma(w)$ or $\sigma(u) < 
\sigma(v) < \sigma(w)$. Now define $\Sigma$ as follows. 
$$
\Sigma_{ij} = \begin{cases}
5 & \mbox{if } i = \sigma(v), j = \sigma(v), \cr 
1 & \mbox{if } i = j, i \neq \sigma(v), \cr 
1 & \mbox{if } i = \sigma(v), j = \sigma(u) \mbox{ or } i = \sigma(v), j = \sigma(w) \cr 
 & \mbox{or } i = \sigma(u), j = \sigma(v) \mbox{ or } i = \sigma(w), j = \sigma(v), \cr 
0 & \mbox{otherwise}. 
\end{cases}
$$

\noindent
Then $\Sigma \in P_{G_\sigma}$. Note that all the diagonal entries of $\Sigma$ are $1$ and 
all off-diagonal entries are $0$ except the $3 \times 3$ submatrix for $\sigma(u), \sigma(v), 
\sigma (w)$. Hence, $\Sigma$ is a permuted block diagonal matrix with $\sigma(u), 
\sigma(v), \sigma(w)$ forming one block and every other index forming a block by itself. 
Using the simple fact that the inverse of a permuted block triangular matrix is permuted 
block triangular, we get that 
$$
\Sigma^{-1}_{ij} = \begin{cases}
\frac{1}{3} & \mbox{if } i = \sigma(v), j = \sigma(v), \cr 
\frac{4}{3} & \mbox{if } i = \sigma(u), j = \sigma(u) \mbox{ or } i = \sigma(w), j = 
\sigma(w), \cr 
1 & \mbox{if } i = j, i \neq \sigma(v) \mbox{ or } \sigma(u) \mbox{ or } \sigma(w), \cr 
-\frac{1}{3} & \mbox{if } i = \sigma(v), j = \sigma(u) \mbox{ or } i = \sigma(v), j = 
\sigma(w), \cr 
 & \mbox{ or } i = \sigma(u), j = \sigma(v) \mbox{ or } i = \sigma(w), j = \sigma(v), \cr 
\frac{1}{3} & \mbox{ if } i = \sigma(u), j = \sigma(w) \mbox{ or } i = \sigma(w), j = 
\sigma(u), \cr 
0 & \mbox{otherwise}. 
\end{cases}
$$

\noindent
Let $C$ denote the maximal clique of $G$ containing $u$ and $v$. Note that $w \notin C$. 
Let $\Sigma_3$ denote the $3 \times 3$ submatrix of $\Sigma$ corresponding to $\sigma(u), 
\sigma(v), \sigma(w)$. Let $\Sigma_3 = L_3 D_3 L_3^T$ denote the modified Cholesky 
decomposition of $\Sigma_3$, and $\Sigma = LDL^T$ be the modified Cholesky decomposition of 
$\Sigma$. For $i,j \in \{\sigma(u), \sigma(v), \sigma(w)\}$, let us define for 
simplicity of notation, $(L_3)_{ij}$ as the entry in the row corresponding to 
$\sigma^{-1} (i)$ and the column corresponding to $\sigma^{-1} (j)$ in $L_3$. Using the 
property that all the diagonal entries of $\Sigma$ are $1$ and all off-diagonal entries are 
$0$ except for $\Sigma_3$, and the uniqueness of the modified Cholesky decomposition of 
$\Sigma$, it follows that 
$$
L_{ij} = \begin{cases}
(L_3)_{ij} & \mbox{if } i > j, \; i,j \in \{\sigma(u), \sigma(v), \sigma(w)\}, \cr 
1 & \mbox{if } i = j, \cr 
0 & \mbox{otherwise}, 
\end{cases}
$$

\noindent
and 
$$
D_{ii} = \begin{cases}
(D_3)_{ii} & \mbox{if } i = \sigma(u), \sigma(v) \mbox{ or } \sigma(w), \cr
1 & \mbox{otherwise}. 
\end{cases}
$$

\noindent
The actual values of the elements of $L_3$ and $D_3$ however, depends on the relative order of $\sigma(u), 
\sigma(v), \sigma(w)$. If $\sigma(v) < \sigma(u) < \sigma(w)$, then $D_{\sigma(v) 
\sigma(v)} = 5, \; D_{\sigma(u) \sigma(u)} = \frac{4}{5}, \; D_{\sigma(w) \sigma(w)} = 
\frac{3}{4}$ and $D_{ii} = 1$ if $i \neq \sigma(v), \sigma(u)$ or $\sigma(w)$. Hence, 
$$
\left| (\Sigma^{-1})_{\sigma(C)} \right| = \frac{1}{3} \neq \prod_{i \in \sigma(C)} 
\frac{1}{D_{ii}} = \frac{1}{4}. 
$$

\noindent
If $\sigma(u) < \sigma(v) < \sigma(w)$, then $D_{\sigma(u) \sigma(u)} = 1, \; D_{\sigma(v) 
\sigma(v)} = 4, \; D_{\sigma(w) \sigma(w)} = \frac{3}{4}$ and $D_{ii} = 1$ if $i \neq 
\sigma(u), \sigma(v)$ or $\sigma(w)$. Hence, 
$$
\left| (\Sigma^{-1})_{\sigma(C)} \right| = \frac{1}{3} \neq \prod_{i \in \sigma(C)}
\frac{1}{D_{ii}} = \frac{1}{4}. 
$$

\noindent
{\bf Case II}: $G$ is homogeneous but $\sigma$ is not a Hasse tree based elimination 
scheme. 

\noindent
Since $\sigma$ is not a Hasse tree based elimination scheme, there exist vertices $a,b \in 
V$ such that $b$ is an ancestor of $a$ in the Hasse tree of $G$, and $\sigma(b) < 
\sigma(a)$. Since $b$ is an ancestor of $a$, there exists $c \in V$, such that $(b,c) 
\in E$ and $(a,c) \notin E$. Now there are three possibilities for the way $\sigma$ orders 
$a,b,c$ given that $\sigma(b) < \sigma(a)$, namely, $\sigma(b) < \sigma(a) < \sigma(c)$ or 
$\sigma(b) < \sigma(c) < \sigma(a)$ or $\sigma(c) < \sigma(b) < \sigma(a)$. Let $v = b, \; u = 
a, \; w = c$ for the first possibility, and $v = b,\; u = c, \; w = a$ for the latter two 
possibilities. Then note that $(u,v), (v,w) \in E, \; (u,w) \notin E$, and $\sigma(v) < 
\sigma(u) < \sigma(w)$ or $\sigma(u) < \sigma(v) < \sigma(w)$. We have thus shown the 
existence of vertices $u,v,w$ such that $(u,v), (v,w) \in E, \; (u,w) \notin E$, and 
$\sigma(v) < \sigma(u) < \sigma(w)$ or $\sigma(u) < \sigma(v) < \sigma(w)$. We can 
therefore use the same $\Sigma$ and maximal clique $C$ as in Case I above, and reach the 
desired conclusion. Hence the result is proved. \hfill$\Box$ 

\vspace{0.1in}

\noindent
We now illustrate the proposition through an example. 
\begin{example}
Consider the homogeneous graph $G$ in Figure \ref{figure2} {\it (a)}. The maximal cliques 
are given by $C_1 = \{w,v',u',u\}$ and $C_2 = \{w,v\}$. The ordering $\sigma$ given by 
$\sigma(w) = 5, \sigma(v) = 4, \sigma(u') = 3, \sigma(u) = 2, \sigma(v') = 1$ is a Hasse 
tree based elimination scheme. Let 
$$
\Sigma = \left( \begin{matrix}
1 & 1 & 1 & 0 & 1 \cr
1 & 2 & 2 & 0 & 2 \cr
1 & 2 & 3 & 0 & 3 \cr
0 & 0 & 0 & 1 & 1 \cr
1 & 2 & 3 & 1 & 5 
\end{matrix} \right) \in P_{G_\sigma}. 
$$

\noindent
Then, 
$$
\left| (\Sigma^{-1})_{\sigma(C_1)} \right| = 1 = \prod_{i \in \sigma(C_1)} \frac{1}{D_{ii}}, 
$$

\noindent
and 
$$
\left| (\Sigma^{-1})_{\sigma(C_2)} \right| = 1 = \prod_{i \in \sigma(C_2)} \frac{1}{D_{ii}}. 
$$

\noindent
Now consider $\sigma$ which is a perfect vertex elimination scheme, but not a Hasse tree 
based elimination scheme, given by $\sigma(v') = 5, \sigma(w) = 4, \sigma(u') = 3, 
\sigma(u) = 2, \sigma(v) = 1$, then 
$$
\Sigma = \left( \begin{matrix}
1 & 0 & 0 & 1 & 0 \cr
0 & 1 & 1 & 1 & 1 \cr
0 & 1 & 2 & 2 & 2 \cr
1 & 1 & 2 & 4 & 3 \cr
0 & 1 & 2 & 3 & 4 
\end{matrix} \right) \in P_{G_\sigma}, 
$$

\noindent
but 
$$
\left| \Sigma^{-1}_{\sigma(C_2)} \right| = 2 \neq \prod_{i \in \sigma(C_2)} 
\frac{1}{D_{ii}} = 1. 
$$

\noindent
Now let $\sigma$ be given by $\sigma(v) = 5, \sigma(u') = 4, \sigma(v') = 3, \sigma(w) = 2, 
\sigma(u) = 1$. Then, $\sigma$ is not a perfect vertex elimination scheme, and 
$$
\Sigma = \left( \begin{matrix}
5 & 1 & 1 & 1 & 0 \cr
1 & 5 & 1 & 1 & 1 \cr
1 & 1 & 5 & 1 & 0 \cr
1 & 1 & 1 & 5 & 0 \cr
0 & 1 & 0 & 0 & 5 
\end{matrix} \right) \in P_{G_\sigma}, 
$$

\noindent
but 
$$
\left| (\Sigma^{-1})_{\sigma(C_1)} \right| = 0.002042484 \neq \prod_{i \in \sigma(C_1)} 
\frac{1}{D_{ii}} = 0.001953125. 
$$

\noindent
Consider the non-homogeneous graph $G$ in Figure \ref{figure2} {\it (b)}. Note however that 
$G$ is a decomposable graph. The maximal cliques are given by $C_1 = \{u',w\}, C_2 = 
\{w,u\}, C_3 = \{u,v\}$. The ordering $\sigma$ given by $\sigma(u') = 4, \sigma(w) = 3, 
\sigma(u) = 2, \sigma(v) = 1$ is a perfect vertex elimination scheme. Let 
$$
\Sigma = \left( \begin{matrix} 
2 & 1 & 0 & 0 \cr 
1 & 2 & 1 & 0 \cr 
0 & 1 & 2 & 1 \cr 
0 & 0 & 1 & 2 
\end{matrix} \right) \in P_{G_\sigma}. 
$$

\noindent
Note however that 
$$
\left| (\Sigma^{-1})_{\sigma(C_3)} \right| = \frac{3}{5} \neq \prod_{i \in \sigma(C_3)} 
\frac{1}{D_{ii}} = \frac{1}{3}. 
$$
\end{example}

\noindent
The two characterizations in the paper are summarized in the main theorem in the introduction. 

\vspace{0.2cm}

\noindent {\it Acknowledgments:} We wish to thank Professor Ingram Olkin for his encouraging remarks on the paper.


\begin{thebibliography}{25}

\bibitem{agler88}
Agler, J., Helton, J.W., McCullough, S. and Rodman, L. (1988). Positive definite matrices with a given sparsity pattern, {\it Linear Algebra Appl.} {\bf 107}, 101--149.

\bibitem{adsnwjrwhc}
Andersson, S. and Wojnar, G.G. (2004). Wishart distributions on homogeneous cones, {\it 
Journal of Theoretical Probability} {\bf 17}, 781-818. 

\bibitem{brndslespn}
Brandstädt, A., Le, V.B. and Spinrad, J. (1999). {\it Graph Classes: A Survey}, SIAM 
Monographs on Discrete Mathematics and Applications. 

\bibitem{dcmplgrchv}
Chvatal, V. (1968). Remark on a paper of Lovasz, {\it Comment. Math. Univ. Carolin.} 
{\bf 9}, 47-50. 

\bibitem{drtonrchsn}
Drton, M. and Richardson, T.S. (2008). Graphical methods for efficient likelihood inference in Gaussian covariance models, {\it Journal of Machine Learning Research} {\bf 9}, 893-914.

\bibitem{dym81}
Dym, H. and Gohberg, I. (1981). Extensions of band matrices with band inverses, {\it Linear Algebra Appl.} {\bf 36}, 1--24.

\bibitem{trnstgphgl}
Gallai, T. (1967). Transitiv orientierbare Graphen, {\it Acta Math. Acad. Sci. Hung.} 
{\bf 18}, 25–66. 

\bibitem{ghlahrilcr}
Ghouila-Houri, A. (1962). Caractérisation des graphes non orientés dont on peut orienter 
les arrêtes de manière à obtenir le graphe $D^*$une relation $D^*$ordre, {\it Les Comptes 
rendus de $L^*$Académie des sciences} {\bf 254}, 1370–1371. 

\bibitem{glmrhfmncg}
Gilmore, P.C. and Hoffman, A.J. (1964). A characterization of comparability graphs and of 
interval graphs, {\it Canadian Journal of Mathematics} {\bf 16}, 539–548. 

\bibitem{glmbcagtpg}
Golumbic, M.C. (1980). {\it Algorithmic Graph Theory and Perfect Graphs}, Academic Press.

\bibitem{grone84}
Gr\H{o}ne, R., Johnson, C.R., Sa, E.M. and Wolkowicz, H. (1984). Positive definite completions of
partial hermitian matrices, {\it Linear Algebra Appl.} {\bf 58}, 109--124.
 
\bibitem{guillot12}
Guillot D., and Rajaratnam, B. (2012). Retaining positive definiteness in thresholded matrices, {\it Linear Algebra Appl.} {\bf 436}, 4143–-4160. 


\bibitem{khrrjtmcvt}
Khare, K. and Rajaratnam, B. (2010). Covariance trees and Wishart distributions on cones, 
{\it Algebraic methods in statistics and probability {II} (AMS CONM SERIES)} {\bf 516}, 
215-223. 

\bibitem{khrrjtmcwp}
Khare, K. and Rajaratnam. B. (2009). Wishart distributions for covariance graph models, 
Technical Report, Stanford University. 

\bibitem{khrrjtmwdg}
Khare, K. and Rajaratnam, B. (2011). Wishart distributions for decomposable covariance 
graph models, {\it Annals of Statistics} {\bf 39}, 514-555. 

\bibitem{lrtzngphmd}
Lauritzen, S.L. (1996). {\it Graphical models}, Oxford University Press Inc., New York. 

\bibitem{ltcmssmwdg}
Letac, G. and Massam, H. (2007). Wishart distributions for decomposable graphs, {\it 
Ann. Statist.} {\bf 35}, 1278-1323. 

\bibitem{trgpapcnhr}
Neher, E. (1999). Transformation groups of the Anderson-Perlman cone, {\it Journal of Lie 
Theory} {\bf 9}, 203-213. 

\bibitem{plsnpwrsth}
Paulsen, V.I., Power, S.C. and Smith, R.R. (1989). Schur products and matrix completions, 
{\it J. Funct. Anal.} {\bf 85}, 151-178. 

\bibitem{pearlwrmth}
Pearl, J. and Wermuth, N. (1994). When can association graphs admit a causal interpretation? In {\it Selecting Models from Data: Artificial Intelligence and Statistics IV, Lecture Notes in Statistics} {\bf 89},  205–214, Springer, New York.

\bibitem{pourahmadi}
Pourahmadi, M. (2007). Cholesky decompositions and estimation of a covariance matrix: orthogonality of variance-correlation parameters, {\it Biometrika} {\bf 94}, 
1006–1013. 

\bibitem{rajaratnam08}
Rajaratnam, B., Massam, H. and Carvalho, C. (2008). Flexible covariance estimation in
graphical models, {\it Annals of Statistics} {\bf 36}, pp. 2818--2849.

\bibitem{roveratocs}
Roverato, A. (2000). Cholesky decomposition of a hyper inverse Wishart matrix, 
{\it Biometrika} {\bf 87}, 99–-112. 

\bibitem{trottercps}
Trotter, W.T. (1992). {\it Combinatorics and Partially Ordered Sets — Dimension Theory}, 
Johns Hopkins University Press. 

\bibitem{urrtprtleg}
Urrutia, J. (1989). {\it Partial orders and Euclidean geometry}, Edited by I. Rival, Kluwer 
Academic Publishers. 


\end{thebibliography}
\end{document}